\documentclass[12pt,a4paper]{article}

\usepackage[utf8]{inputenc}
\usepackage[T1]{fontenc}
\usepackage{amsmath,amssymb,amsfonts}
\usepackage{graphicx}
\usepackage{tikz}
\usetikzlibrary{arrows.meta,calc,positioning,decorations.pathreplacing,shapes.misc}
\usepackage{listings}
\usepackage{xcolor}
\usepackage{hyperref}
\usepackage{geometry}
\usepackage{booktabs}
\usepackage{caption}
\usepackage{bm}
\usepackage{float}

\definecolor{matlabbg}{RGB}{248,248,248}

\lstdefinestyle{matlab}{
  language=Matlab,
  basicstyle=\small\ttfamily,
  backgroundcolor=\color{matlabbg},
  keywordstyle=\bfseries,
  commentstyle=\itshape,
  numbers=left,
  numberstyle=\tiny\color{gray},
  numbersep=6pt,
  frame=single,
  breaklines=true,
  breakatwhitespace=false,
  tabsize=4,
  showstringspaces=false,
  captionpos=b,
  xleftmargin=1.5em,
  framexleftmargin=1em,
}
\title{\textbf{Solving Maxwell's Equations with Mimetic Methods}}
\author{
  Johnny Corbino\thanks{Lawrence Berkeley National Laboratory. E-mail: \texttt{jcorbino@lbl.gov}}
}
\date{}

\begin{document}
\maketitle

\begin{abstract}
We present a mimetic finite-difference approach for solving Maxwell's equations in one and two spatial dimensions. After introducing the governing equations and the classical Finite-Difference Time-Domain (FDTD) method, we describe mimetic operators that satisfy a discrete analogue of the extended Gauss divergence theorem and show how they lead to a compact, physically consistent formulation for computational electromagnetics. Two numerical examples are presented: a one-dimensional sinusoidal wave interacting with a lossy dielectric slab, and a two-dimensional Gaussian pulse with Uniaxial Perfectly Matched Layer (UPML) absorbing boundary conditions. All implementations use the Mimetic Operators Library Enhanced (MOLE).
\end{abstract}

% ══════════════════════════════════════════════════════════════
\section{Introduction}\label{sec:intro}

Maxwell's equations govern the behavior of electromagnetic fields and are central to a vast range of engineering and scientific applications. The Finite-Difference Time-Domain (FDTD) method, introduced by Yee~\cite{Yee1966}, remains one of the most popular techniques for the numerical solution of these equations due to its explicit time-marching nature and conceptual simplicity~\cite{Taflove2005,Sullivan2013}. However, the standard FDTD scheme is limited to second-order spatial accuracy, and extending it to higher orders while preserving the underlying physics requires considerable effort.

Mimetic finite-difference methods offer an alternative that addresses these limitations. They construct discrete operators---divergence, gradient, curl, and Laplacian---that satisfy the same vector calculus identities as their continuous counterparts, including a discrete version of the extended Gauss divergence theorem~\cite{CastilloGrone2003,Castillo2013,Corbino2020}. This ensures conservation properties at the discrete level and provides uniform accuracy up to and including the domain boundaries.

This paper is organized as follows. Section~\ref{sec:maxwell} introduces Maxwell's equations and defines all relevant physical parameters. Section~\ref{sec:fdtd} gives a brief overview of the FDTD method. Section~\ref{sec:mimetic} presents the mimetic operator framework. Sections~\ref{sec:example1d} and~\ref{sec:example2d} present one-dimensional and two-dimensional numerical examples, respectively. Finally, Section~\ref{sec:conclusions} offers concluding remarks.

% ══════════════════════════════════════════════════════════════
\section{Maxwell's Equations}\label{sec:maxwell}

In differential form, for linear, isotropic media, Maxwell's equations read:
\begin{align}
  \nabla \times \mathbf{E} &= -\mu \frac{\partial \mathbf{H}}{\partial t}, \label{eq:faraday}\\[6pt]
  \nabla \times \mathbf{H} &= \varepsilon \frac{\partial \mathbf{E}}{\partial t} + \sigma \mathbf{E}, \label{eq:ampere}\\[6pt]
  \nabla \cdot (\varepsilon \mathbf{E}) &= \rho_v, \label{eq:gauss_e}\\[6pt]
  \nabla \cdot (\mu \mathbf{H}) &= 0, \label{eq:gauss_m}
\end{align}
where:
\begin{itemize}
  \item $\mathbf{E}$ is the electric field intensity (V/m),
  \item $\mathbf{H}$ is the magnetic field intensity (A/m),
  \item $\varepsilon$ is the electric permittivity (F/m); in free space $\varepsilon_0 \approx 8.854 \times 10^{-12}$~F/m, and for a material with relative permittivity $\varepsilon_r$ we have $\varepsilon = \varepsilon_r\,\varepsilon_0$,
  \item $\mu$ is the magnetic permeability (H/m); in free space $\mu_0 = 4\pi \times 10^{-7}$~H/m, and for a material with relative permeability $\mu_r$ we have $\mu = \mu_r\,\mu_0$,
  \item $\sigma$ is the electric conductivity (S/m), modeling ohmic losses ($\sigma = 0$ for lossless media),
  \item $\rho_v$ is the volume charge density (C/m$^3$).
\end{itemize}

The phase velocity in a medium is $v = 1/\sqrt{\mu\varepsilon}$; in free space, $v = c_0 = 1/\sqrt{\mu_0\varepsilon_0} \approx 3\times10^8$~m/s (speed of light in vacuum).

In one dimension, restricting to waves propagating along $z$ with $E_x$ and $H_y$, Faraday's and Amp\`ere's laws reduce to:
\begin{align}
  \varepsilon \frac{\partial E_x}{\partial t} &= -\frac{\partial H_y}{\partial z} - \sigma E_x, \label{eq:1d_E}\\[4pt]
  \mu\frac{\partial H_y}{\partial t} &= -\frac{\partial E_x}{\partial z}. \label{eq:1d_H}
\end{align}

The magnetic flux density $\mathbf{B}$ and the magnetic field intensity $\mathbf{H}$ are related by
\begin{equation}\label{eq:B_H}
  \mathbf{B} = \mu\,\mathbf{H}.
\end{equation}
In the 1D example that follows, we work with $\mathbf{H}$, which is the natural variable in Amp\`ere's law. In the 2D example, we work with $\mathbf{B}$, which simplifies the formulation when $\mu$ is constant (free space, $\mu = \mu_0$) because Faraday's law becomes $\partial\mathbf{B}/\partial t = -\nabla\times\mathbf{E}$ without a material coefficient.

For two-dimensional problems in the TM$_z$ polarization, the relevant fields are $E_z$, $B_x$, and $B_y$, and Faraday's and Amp\`ere's laws reduce to:
\begin{align}
  \frac{\partial B_x}{\partial t} &= -\frac{\partial E_z}{\partial y}, \label{eq:2d_Bx}\\[4pt]
  \frac{\partial B_y}{\partial t} &= \phantom{-}\frac{\partial E_z}{\partial x}, \label{eq:2d_By}\\[4pt]
  \varepsilon\frac{\partial E_z}{\partial t} &= \frac{1}{\mu}\!\left(\frac{\partial B_y}{\partial x} - \frac{\partial B_x}{\partial y}\right). \label{eq:2d_Ez}
\end{align}

% ══════════════════════════════════════════════════════════════
\section{The FDTD Method}\label{sec:fdtd}

The FDTD method discretizes Maxwell's equations on a \emph{staggered grid} introduced by Yee~\cite{Yee1966}. Electric and magnetic field components are sampled at spatially offset locations, and their time evolution is staggered by half a time step, yielding the \emph{leapfrog} scheme.

\begin{figure}[H]
\centering
\begin{tikzpicture}[>=Stealth, scale=1.0]
  % E-field nodes
  \foreach \i in {0,1,2,3,4,5} {
    \filldraw (\i*1.8, 0) circle (2.5pt);
    \node[below=6pt] at (\i*1.8, 0) {\small $E_x^n(\i)$};
  }
  % H-field nodes
  \foreach \i in {0,1,2,3,4} {
    \draw (\i*1.8+0.9, 0.6) node[cross out, draw, thick, inner sep=2.5pt] {};
    \node[above=6pt] at (\i*1.8+0.9, 0.6) {\small $H_y^{n+\frac{1}{2}}$};
  }
  % Grid line
  \draw[gray, thin] (0,0) -- (9,0);
  % Spacing
  \draw[<->] (0, -0.9) -- (1.8, -0.9);
  \node[below] at (0.9, -0.9) {$\Delta z$};
  \draw[<->] (0, 0.15) -- (0.9, 0.15);
  \node[above] at (0.45, 0.15) {\small $\frac{\Delta z}{2}$};
\end{tikzpicture}
\caption{One-dimensional Yee staggered grid. Electric field values (dots) are stored at integer nodes; magnetic field values (crosses) at half-integer nodes.}
\label{fig:yee_1d}
\end{figure}
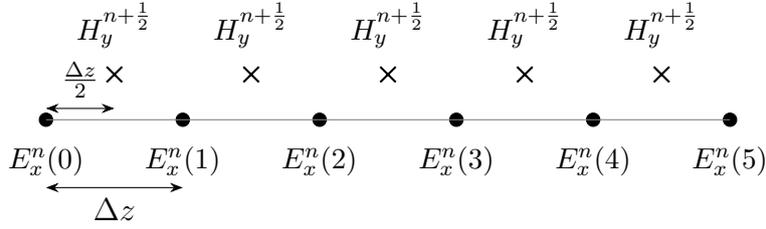

For the one-dimensional case, in normalized units ($c_0=1$, $\mu_0=1$, $\varepsilon_0=1$, so that $\varepsilon = \varepsilon_r$ and $\mu = \mu_r$), the update equations are:
\begin{align}
  E_x^{n+1}(i) &= E_x^{n}(i) + \frac{\Delta t}{\varepsilon\,\Delta z}\left[H_y^{n+1/2}\!\left(i\!-\!\tfrac{1}{2}\right) - H_y^{n+1/2}\!\left(i\!+\!\tfrac{1}{2}\right)\right], \label{eq:fdtd_E}\\[6pt]
  H_y^{n+3/2}\!\left(i\!+\!\tfrac{1}{2}\right) &= H_y^{n+1/2}\!\left(i\!+\!\tfrac{1}{2}\right) + \frac{\Delta t}{\mu\,\Delta z}\left[E_x^{n+1}(i) - E_x^{n+1}(i\!+\!1)\right]. \label{eq:fdtd_H}
\end{align}

The scheme is stable under the Courant--Friedrichs--Lewy (CFL) condition $\Delta t \leq \Delta z / c_0$.

In a lossy dielectric with relative permittivity $\varepsilon_r$ and conductivity $\sigma$, the electric field update acquires damping coefficients. Using a semi-implicit discretization of the conduction current term, one obtains:
\begin{equation}\label{eq:fdtd_E_lossy}
  E_x^{n+1}(i) = c_a(i)\,E_x^{n}(i) + c_b(i)\left[H_y^{n+1/2}\!\left(i\!-\!\tfrac{1}{2}\right) - H_y^{n+1/2}\!\left(i\!+\!\tfrac{1}{2}\right)\right],
\end{equation}
where
\begin{equation}\label{eq:ca_cb}
  c_a = \frac{1 - \dfrac{\sigma\,\Delta t}{2\,\varepsilon}}{1 + \dfrac{\sigma\,\Delta t}{2\,\varepsilon}}, \qquad
  c_b = \frac{\dfrac{\Delta t}{\varepsilon\,\Delta z}}{1 + \dfrac{\sigma\,\Delta t}{2\,\varepsilon}}.
\end{equation}
The coefficient $c_a < 1$ causes exponential decay of the field inside the lossy region, while $c_b$ controls the strength of the electric field response to spatial variations in the magnetic field.

While FDTD is straightforward and efficient, it has important limitations. The Yee scheme does preserve the discrete divergence conditions ($\nabla\cdot(\varepsilon\mathbf{E})=\rho_v$ and $\nabla\cdot\mathbf{B}=0$) if they are satisfied initially, thanks to the staggered grid structure. However, when extending to higher orders of accuracy or different coordinate systems (e.g., curvilinear), the extension is not trivial and often a divergence correction method is required~\cite{Munz2000}.

% ══════════════════════════════════════════════════════════════
\section{Mimetic Finite-Difference Operators}\label{sec:mimetic}

\subsection{Definitions and properties}

Mimetic operators construct discrete analogues of the differential operators---divergence ($\mathbf{D} \equiv \nabla\cdot$), gradient ($\mathbf{G} \equiv \nabla$), curl ($\mathbf{C} \equiv \nabla\times$), and Laplacian ($\mathbf{L} \equiv \nabla^2$)---that satisfy the fundamental vector calculus identities~\cite{Castillo2013,Corbino2020}. These identities are preserved at any order of accuracy and extend naturally to curvilinear coordinate systems~\cite{Boada2025}, addressing the limitations of the Yee scheme discussed above:
\begin{align}
  \mathbf{G}\,f_{\text{const}} &= \mathbf{0}, \label{eq:mim_id1}\\
  \mathbf{D}\,\boldsymbol{\nu}_{\text{const}} &= \mathbf{0}, \label{eq:mim_id2}\\
  \mathbf{C}\,\mathbf{G}\,f &= \mathbf{0}, \label{eq:mim_id3}\\
  \mathbf{D}\,\mathbf{C}\,\boldsymbol{\nu} &= \mathbf{0}, \label{eq:mim_id4}\\
  \mathbf{D}\,\mathbf{G}\,f &= \mathbf{L}\,f. \label{eq:mim_id5}
\end{align}

Moreover, the Corbino--Castillo (CC) operators satisfy the \emph{discrete extended Gauss divergence theorem}~\cite{Corbino2020}:
\begin{equation}\label{eq:gauss_discrete}
  \langle \mathbf{D}\boldsymbol{\nu}, f \rangle_Q + \langle \mathbf{G}f, \boldsymbol{\nu} \rangle_P = \langle \boldsymbol{\mathcal{B}}\boldsymbol{\nu}, f \rangle,
\end{equation}
where $P$ and $Q$ are diagonal weight matrices and $\boldsymbol{\mathcal{B}}$ is the mimetic boundary operator. This identity is the discrete counterpart of
\begin{equation}
  \int_\Omega (\nabla\cdot\boldsymbol{\nu})\,f\,dV + \int_\Omega \nabla f \cdot \boldsymbol{\nu}\,dV = \oint_{\partial\Omega} f\,\boldsymbol{\nu}\cdot\hat{\mathbf{n}}\,dS,
\end{equation}
and guarantees that the discrete scheme inherits the conservation properties of the continuous equations.

\subsection{Staggered grids}

Like the Yee grid, mimetic operators are defined on staggered grids. In one dimension, scalar quantities are stored at cell centers and at the two boundary nodes, while vector components reside at cell edges (nodes).

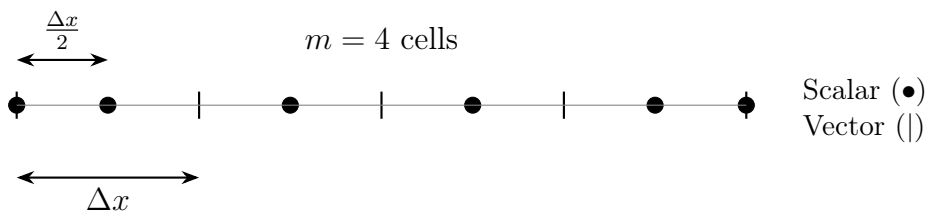
\begin{figure}[H]
\centering
\begin{tikzpicture}[>=Stealth, scale=1.2]
  % Cell edges (nodes for H)
  \foreach \i in {0,1,2,3,4} {
    \draw[thick] (\i*2, -0.15) -- (\i*2, 0.15);
  }
  % Cell centers + boundary nodes (for E)
  \filldraw (0, 0) circle (2.5pt);
  \filldraw (8, 0) circle (2.5pt);
  \foreach \i in {0,1,2,3} {
    \filldraw (\i*2+1, 0) circle (2.5pt);
  }
  % Grid line
  \draw[gray] (0,0) -- (8,0);
  % Spacing
  \draw[<->, thick] (0, -0.8) -- (2, -0.8);
  \node[below] at (1, -0.8) {$\Delta x$};
  \draw[<->, thick] (0, 0.5) -- (1, 0.5);
  \node[above] at (0.5, 0.5) {$\frac{\Delta x}{2}$};
  % Legend
  \node[right] at (8.5, 0.15) {\small Scalar ($\bullet$)};
  \node[right] at (8.5, -0.25) {\small Vector ($|$)};
  \node[above=12pt] at (4, 0.15) {$m = 4$ cells};
\end{tikzpicture}
\caption{One-dimensional mimetic staggered grid with $m=4$ cells. Scalar values (dots) reside at cell centers and boundary nodes ($m+2$ locations). Vector values (ticks) reside at cell edges ($m+1$ locations).}
\label{fig:mimetic_1d}
\end{figure}

For $m$ cells, the \textbf{gradient} $\mathbf{G}$ is an $(m+1)\times(m+2)$ sparse matrix and the \textbf{divergence} $\mathbf{D}$ is an $(m+2)\times(m+1)$ sparse matrix.

\subsection{One-dimensional operators}

The second-order mimetic gradient, following Corbino and Castillo~\cite{Corbino2020}, is:
\begin{equation}\label{eq:G2}
  \mathbf{G} = \frac{1}{\Delta x}
  \begin{bmatrix}
    -\frac{8}{3} & 3 & -\frac{1}{3} \\
     & -1 & 1 \\
     & & \ddots & \ddots \\
     & & & -1 & 1 \\
     & & & \frac{1}{3} & -3 & \frac{8}{3}
  \end{bmatrix}_{(m+1)\times(m+2)}.
\end{equation}

The mimetic operators are constructed via the solution of Vandermonde systems and achieve a uniform order of accuracy throughout the domain, including at the boundaries. This distinguishes them from standard finite-difference and summation-by-parts (SBP) operators, which typically suffer from reduced accuracy near the boundaries.

The second-order mimetic divergence is:
\begin{equation}\label{eq:D2}
  \mathbf{D} = \frac{1}{\Delta x}
  \begin{bmatrix}
    0 & \cdots & 0 \\
    -1 & 1 & \\
     & \ddots & \ddots \\
     & -1 & 1 \\
    0 & \cdots & 0
  \end{bmatrix}_{(m+2)\times(m+1)}.
\end{equation}

\subsection{Two-dimensional operators via Kronecker products}

Two-dimensional operators are constructed from one-dimensional ones using Kronecker products:
\begin{equation}\label{eq:G2D}
  \mathbf{G}_{xy} = \begin{bmatrix} \hat{I}_n^T \otimes \mathbf{G}_x \\ \mathbf{G}_y \otimes \hat{I}_m^T \end{bmatrix}, \qquad
  \mathbf{D}_{xy} = \begin{bmatrix} \hat{I}_n \otimes \mathbf{D}_x & \mathbf{D}_y \otimes \hat{I}_m \end{bmatrix},
\end{equation}
where $\hat{I}_m$ is an $(m+2)\times m$ augmented identity (zeros in the first and last rows).

\begin{figure}[H]
\centering
\begin{tikzpicture}[>=Stealth, scale=0.85]
  % Draw cells
  \foreach \i in {0,1,2,3} {
    \foreach \j in {0,1,2} {
      \draw[gray, thin] (\i*2, \j*2) rectangle (\i*2+2, \j*2+2);
    }
  }
  % Scalar nodes at cell centers
  \foreach \i in {0,1,2,3} {
    \foreach \j in {0,1,2} {
      \filldraw (\i*2+1, \j*2+1) circle (2.5pt);
    }
  }
  % Boundary scalar nodes
  \foreach \j in {0,1,2} {
    \filldraw (0, \j*2+1) circle (2.5pt);
    \filldraw (8, \j*2+1) circle (2.5pt);
  }
  \foreach \i in {0,1,2,3} {
    \filldraw (\i*2+1, 0) circle (2.5pt);
    \filldraw (\i*2+1, 6) circle (2.5pt);
  }
  \filldraw (0,0) circle (2.5pt);
  \filldraw (8,0) circle (2.5pt);
  \filldraw (0,6) circle (2.5pt);
  \filldraw (8,6) circle (2.5pt);
  % Vector components on edges (short lines)
  \foreach \i in {0,1,2,3} {
    \foreach \j in {0,1,2,3} {
      \draw[thick] (\i*2+1, \j*2-0.3) -- (\i*2+1, \j*2+0.3);
    }
  }
  \foreach \i in {0,1,2,3,4} {
    \foreach \j in {0,1,2} {
      \draw[thick] (\i*2-0.3, \j*2+1) -- (\i*2+0.3, \j*2+1);
    }
  }
  % Labels
  \draw[<->] (0, -0.5) -- (2, -0.5);
  \node[below] at (1, -0.5) {$\Delta x$};
  \draw[<->] (-0.5, 0) -- (-0.5, 2);
  \node[left] at (-0.5, 1) {$\Delta y$};
\end{tikzpicture}
\caption{Two-dimensional staggered grid ($m=4$, $n=3$). Scalar values (dots) at cell centers and boundaries; vector components (ticks perpendicular to cell sides) at cell edges.}
\label{fig:mimetic_2d}
\end{figure}

\subsection{Application to Maxwell's equations}

Using the same normalized units ($c_0=1$, $\mu_0=1$, $\varepsilon_0=1$), the leapfrog update for Maxwell's equations becomes:
\begin{align}
  \mathbf{E}^{n+1} &= \mathbf{E}^{n} - \Delta t\;\frac{1}{\varepsilon_r}\;\mathbf{D}\,\mathbf{H}^{n+1/2}, \label{eq:mim_update_E}\\[4pt]
  \mathbf{H}^{n+3/2} &= \mathbf{H}^{n+1/2} - \Delta t\;\frac{1}{\mu_r}\;\mathbf{G}\,\mathbf{E}^{n+1}. \label{eq:mim_update_H}
\end{align}

Each time step reduces to two sparse matrix--vector products. The operators $\mathbf{D}$ and $\mathbf{G}$ are constructed once, and changing the order of accuracy requires only changing the parameter $k$ passed to the MOLE library~\cite{Corbino2024}. Material properties enter through element-wise scaling of the operator (for permittivity) or multiplicative coefficients (for conductivity), without modifying the operator construction itself.

% ══════════════════════════════════════════════════════════════
\section{1D Example: Sinusoidal Wave in a Lossy Dielectric}\label{sec:example1d}

This example, adapted from Sullivan~\cite{Sullivan2013}, considers a one-dimensional domain with $m = 200$ cells. A continuous sinusoidal source at frequency $f = 700$~MHz is placed near the left boundary. A lossy dielectric slab with relative permittivity $\varepsilon_r = 4$ and conductivity $\sigma = 0.04$~S/m occupies the right half of the domain (from cell 100 onward). Absorbing boundary conditions are applied at both ends.

\subsection{Loss coefficients}

Following the semi-implicit discretization of the conduction current (see Section~\ref{sec:fdtd}), the loss coefficients are:
\begin{equation}
  c_a = \frac{1 - \ell}{1 + \ell}, \qquad
  c_b = \frac{1/(2\varepsilon_r)}{1 + \ell}, \qquad
  \ell = \frac{\sigma\,\Delta t}{2\,\varepsilon_0\,\varepsilon_r},
\end{equation}
where $\Delta t = 0.01/(2c_0)$ is the time step, chosen so that the Courant number is $S_c = c_0\,\Delta t/\Delta z = 0.5$. In free space, $c_a = 1$ and $c_b = 0.5$.

\subsection{Mimetic implementation}

The mimetic operators $\mathbf{D}$ and $\mathbf{G}$ are built with unit spacing; the loss coefficients $c_a$ and $c_b$ encode both permittivity and conductivity. A simple absorbing boundary condition (ABC) is used at both ends: it saves the adjacent interior values before the update and overwrites the boundary entries afterward.
\begin{lstlisting}[caption={1D mimetic Maxwell solver: sinusoidal wave in a lossy dielectric.}]
m = 200; % Number of cells
k = 2; % Order of accuracy

% Mimetic operators
D = div(k, m, 1);
G = grad(k, m, 1);

dt    = 0.01/(2*3e8); % Time step
steps = 500;

% Initialize fields
ex = zeros(m+2, 1);
hy = zeros(m+1, 1);

% Wave parameters
freq = 700e6;

% Lossy dielectric slab
kstart  = 100;
epsilon = 4;
sigma   = 0.04;
epsz    = 8.85419e-12;
loss_term = dt*sigma/(2*epsz*epsilon);

cb = 0.5*ones(m+2, 1);
ca = ones(m+2, 1);
ca(kstart:end) = (1 - loss_term)/(1 + loss_term);
cb(kstart:end) = (0.5/epsilon)/(1 + loss_term);

% Time loop
for n = 1:steps
    El = ex(2);       % Save boundary values
    Er = ex(end-1);

    ex = ca.*ex - cb.*(D*hy);              % E update
    ex(5) = ex(5) + sin(2*pi*freq*dt*n);   % Source
    ex(1) = El;  ex(end) = Er;             % ABC

    hy = hy - 0.5*(G*ex);                  % H update
end
\end{lstlisting}

The update \texttt{ex = ca.*ex - cb.*(D*hy)} is the mimetic analogue of~\eqref{eq:fdtd_E_lossy}. The element-wise products with $c_a$ and $c_b$ apply the material properties point-wise, while the matrix--vector product $\mathbf{D}\,\mathbf{H}$ computes the spatial derivative via the mimetic divergence operator. Inside the lossy slab, $c_a < 1$ damps the existing field at each time step, while $c_b < 0.5$ reduces the wave speed, shortening the wavelength by a factor of $\sqrt{\varepsilon_r}$. Together, the two coefficients produce the exponential attenuation characteristic of lossy media.  

\begin{figure}[H]
\centering
\includegraphics[width=0.85\textwidth]{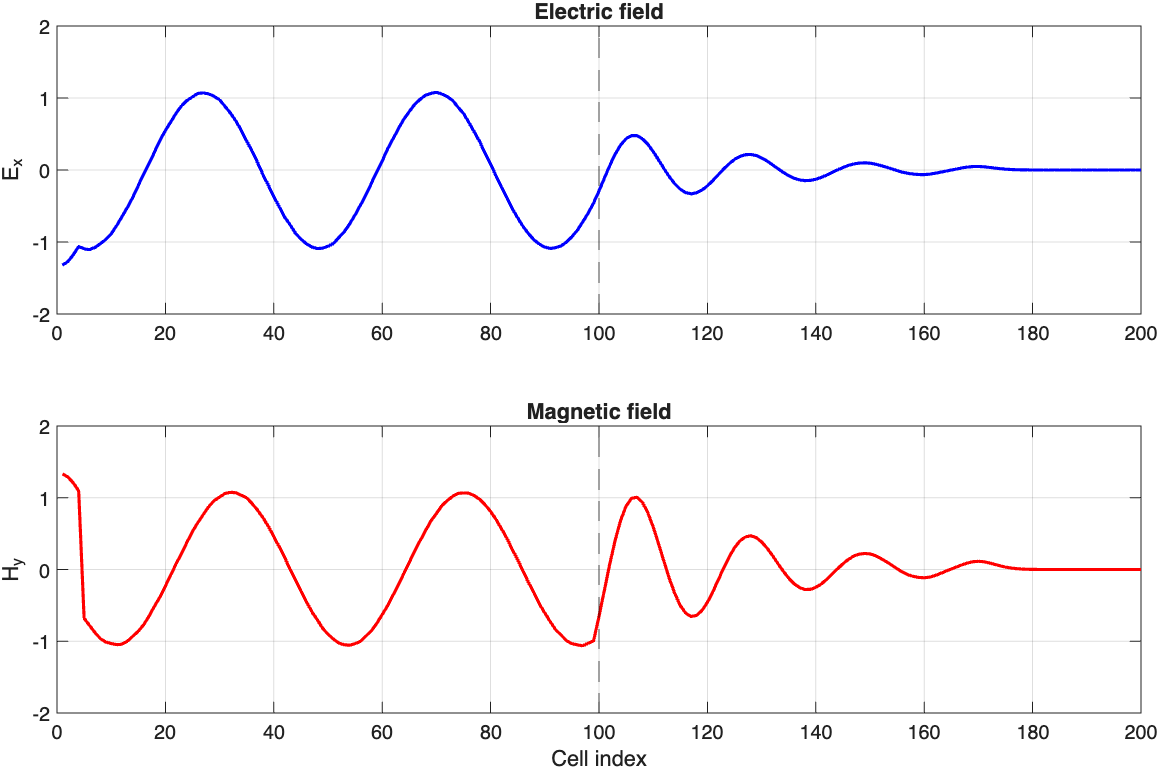}
\caption{Electric and magnetic fields after 500 steps for a 700~MHz sinusoidal wave interacting with a lossy dielectric slab ($\varepsilon_r = 4$, $\sigma = 0.04$~S/m).}
\label{fig:1d_result}
\end{figure}

% ══════════════════════════════════════════════════════════════
\section{2D Example: Gaussian Pulse with UPML}\label{sec:example2d}

\subsection{Problem setup}

This example is also adapted from Sullivan~\cite{Sullivan2013}. We consider a $100 \times 100$ cell domain $[0,1]^2$ with $\Delta x = \Delta y = 0.01$ and $\Delta t = 0.5\,\min(\Delta x,\Delta y)$. A Gaussian pulse $E_z(x,y,0) = \exp\!\bigl[-400\bigl((x-0.5)^2+(y-0.5)^2\bigr)\bigr]$ is placed at the center and radiates outward.

\subsection{Uniaxial Perfectly Matched Layer (UPML)}

The UPML surrounds the computational domain with a layer of $d_{\text{PML}}$ cells in which the artificial conductivities $\sigma_x(x)$ and $\sigma_y(y)$ ramp up smoothly:
\begin{equation}\label{eq:sigma_profile}
  \sigma_x(i) = \sigma_{\max}\!\left(\frac{\text{depth}}{d_{\text{PML}}}\right)^{\!p},
\end{equation}
where $\text{depth}$ is the distance into the PML in grid cells, $\sigma_{\max}$ is the maximum conductivity, and $p$ is a polynomial grading exponent that controls the smoothness of the ramp.

The layer is \emph{perfectly matched} in the sense that the modified medium preserves the intrinsic impedance $\eta = \sqrt{\mu/\varepsilon}$ of the original domain. Because $\eta$ is continuous across the interface, an incoming plane wave sees no impedance mismatch and enters the absorbing layer without reflection, regardless of its polarization, speed, or angle of incidence. Attenuation is introduced solely through the artificial conductivity~$\sigma$, which damps the wave amplitude without altering the wave impedance. A common estimate for the maximum conductivity is~\cite{Taflove2005}
\begin{equation}\label{eq:sigma_opt}
  \sigma_{\max} = \frac{0.8\,(p+1)}{\eta\,\Delta x}.
\end{equation}
In normalized units ($\mu_0 = \varepsilon_0 = 1$), the free-space impedance $\eta_0 = \sqrt{\mu_0/\varepsilon_0} = 1$ and $\sigma_{\max}$ depends only on the grading exponent and grid spacing. The value $\sigma_{\max} = 100$ used here is a practical choice that provides sufficient absorption while avoiding excessive damping near the PML interface.

The total 2D conductivity is formed as a separable sum $\sigma(i,j) = \sigma_x(i) + \sigma_y(j)$, and the damping factors are computed as:
\begin{equation}\label{eq:damping}
  a_E(i,j) = \exp\!\bigl(-\sigma(i,j)\,\Delta t\bigr), \qquad
  a_B = \exp\!\bigl(-\sigma_B\,\Delta t\bigr),
\end{equation}
where $\sigma_B$ is evaluated at the grid locations that correspond to the magnetic field. The UPML-augmented leapfrog update is:
\begin{align}
  \mathbf{B}^{n+1/2} &= a_B \odot \bigl(\mathbf{B}^{n-1/2} - \Delta t\,\mathbf{G}_{2D}\,\mathbf{E}^{n}\bigr), \label{eq:upml_B}\\
  \mathbf{E}^{n+1} &= a_E \odot \bigl(\mathbf{E}^{n} - \Delta t\,\mathbf{D}_{2D}\,\mathbf{B}^{n+1/2}\bigr), \label{eq:upml_E}
\end{align}
where $\odot$ denotes the Hadamard product. The mimetic operators themselves are \emph{unmodified}; the PML enters only through the multiplicative damping vectors.

\subsection{Mimetic implementation}

\begin{lstlisting}[caption={2D mimetic Maxwell solver with UPML (Gaussian pulse).}]
mx = 100; my = 100; % Number of cells
dx = 1/mx; dy = 1/my; % Grid spacing
dt = 0.5*min(dx,dy); % Time step
steps = 140;

% Mimetic operators
k = 2;
G = dt*grad2D(k, mx, dx, my, dy);
D = dt*div2D(k, mx, dx, my, dy);

% Gaussian initial condition
E = exp(-400*((XE-0.5).^2 + (YE-0.5).^2));
E = E(:);
B = zeros(NBx+NBy, 1);

% UPML profiles (polynomial grading)
pml = 30;  sigma_max = 100;  p = 4;
sigma_x = zeros(mx+2,1);
sigma_y = zeros(my+2,1);

Lx = ix <= pml;   Rx = ix >= mx+3-pml;
Ly = iy <= pml;   Ry = iy >= my+3-pml;

sigma_x(Lx) = sigma_max*((pml-ix(Lx)+1)/pml).^p;
sigma_x(Rx) = sigma_max*((ix(Rx)-(mx+2-pml))/pml).^p;
sigma_y(Ly) = sigma_max*((pml-iy(Ly)+1)/pml).^p;
sigma_y(Ry) = sigma_max*((iy(Ry)-(my+2-pml))/pml).^p;

% 2D damping (separable outer sum)
SIGE = sigma_y + sigma_x.';
aE = exp(-SIGE(:)*dt);

sigmaBx = sigma_y(1:my+1) + sigma_x(2:mx+1).';
sigmaBy = sigma_y(2:my+1) + sigma_x(1:mx+1).';
aB = exp(-[sigmaBx(:); sigmaBy(:)]*dt);

% Leapfrog half-step
B = aB.*(B - 0.5*G*E);

% Time loop
for n = 1:steps
    % Update fields
    B = aB.*(B - G*E);
    E = aE.*(E - D*B);
end
\end{lstlisting}

The time loop consists of just two lines: one for the magnetic update and one for the electric update. The damping vectors $a_B$ and $a_E$ smoothly attenuate outgoing waves inside the PML, while the interior of the domain remains unaffected ($a = 1$). The clean separation between operator construction and boundary treatment is a hallmark of the mimetic approach.

\begin{figure}[H]
\centering
\includegraphics[width=0.85\textwidth]{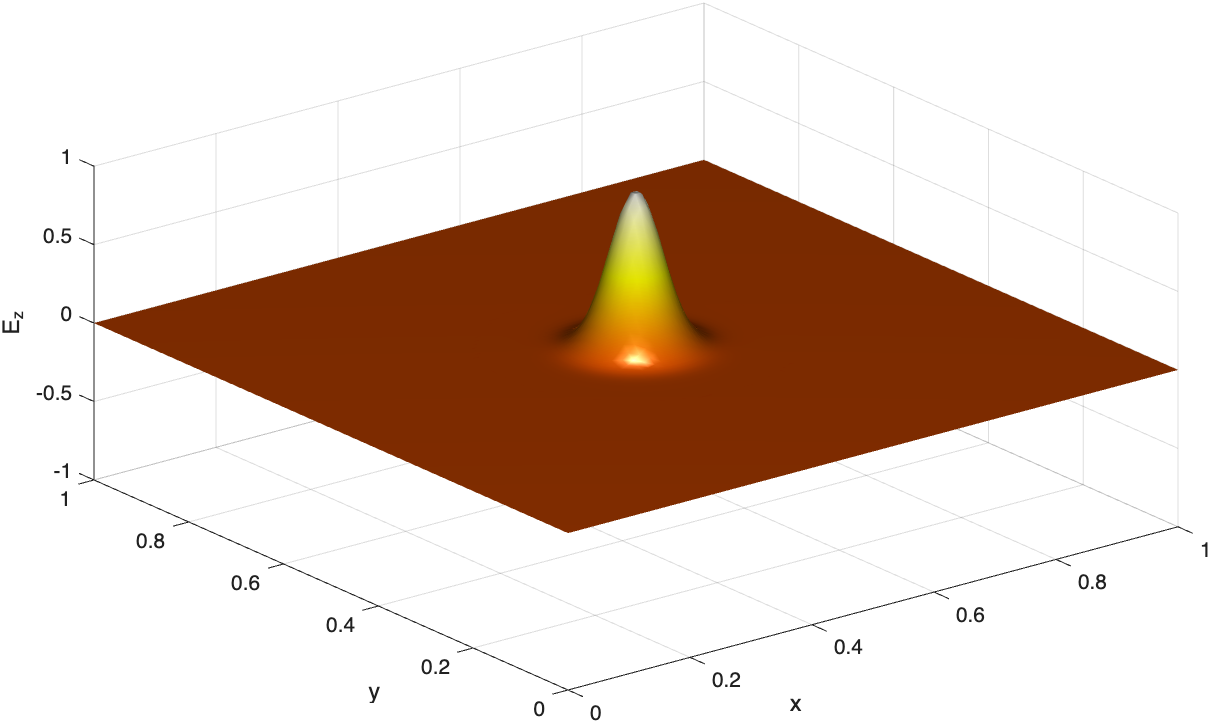}
\caption{Initial electric field ($E_z$) for a 2D Gaussian pulse with UPML absorbing boundaries ($d_{\text{PML}} = 30$ cells, $\sigma_{\max} = 100$, $p = 4$).}
\label{fig:2d_result}
\end{figure}

\begin{figure}[H]
\centering
\includegraphics[width=0.85\textwidth]{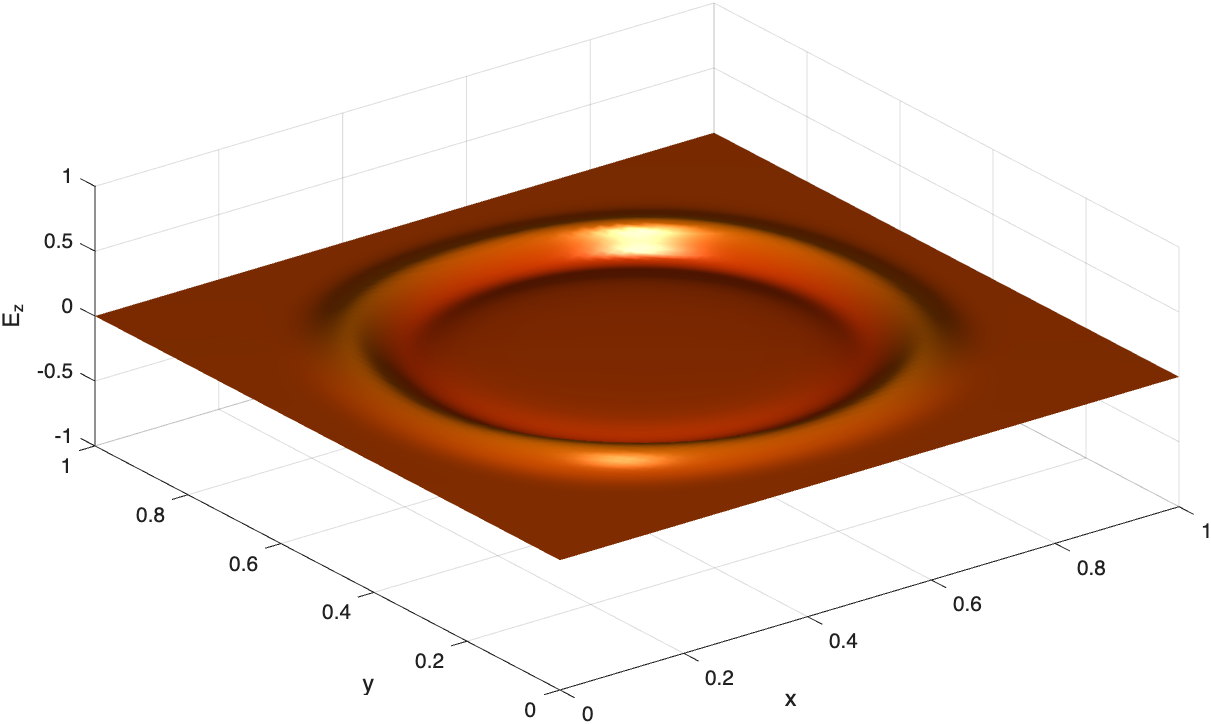}
\caption{Electric field ($E_z$) after 70 steps.}
\label{fig:2d_result_2}
\end{figure}

\begin{figure}[H]
\centering
\includegraphics[width=0.85\textwidth]{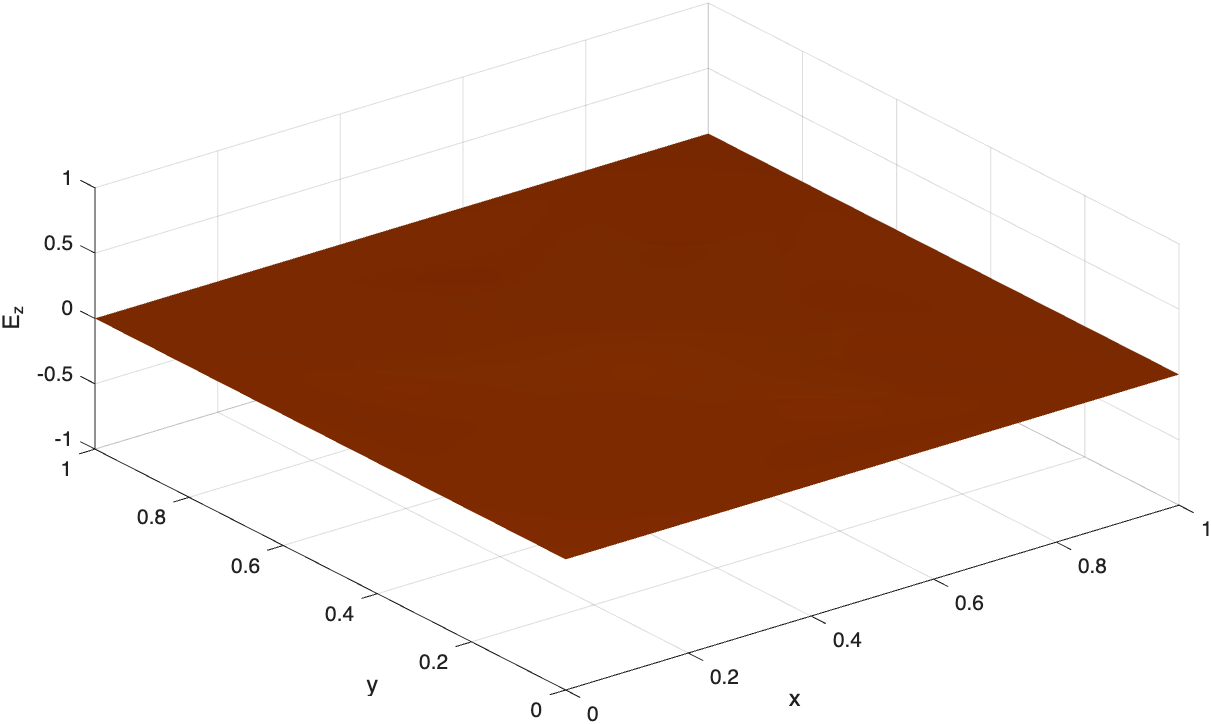}
\caption{Electric field ($E_z$) after 140 steps.}
\label{fig:2d_result_3}
\end{figure}

% ══════════════════════════════════════════════════════════════
\section{Concluding Remarks}\label{sec:conclusions}

We have presented a mimetic finite-difference framework for solving Maxwell's equations, demonstrating its application to one- and two-dimensional electromagnetic wave propagation problems. Compared to the classical FDTD method, the mimetic approach offers several advantages:

\begin{enumerate}
  \item \textbf{Operator-based formulation.} The mimetic approach replaces index-level update loops with sparse matrix--vector products. Each time step reduces to two multiplications regardless of dimensionality or order of accuracy, yielding a code that is shorter and less error-prone.

  \item \textbf{Uniform boundary accuracy.} The Corbino--Castillo operators achieve the same order of accuracy at the domain boundaries as in the interior, unlike standard finite-difference or SBP operators. This is critical for wave problems where boundary errors can pollute the entire solution.

  \item \textbf{Discrete conservation.} Fulfillment of the discrete extended Gauss divergence theorem~\eqref{eq:gauss_discrete} ensures that fundamental conservation properties are preserved, leading to better long-term stability and physical fidelity.

  \item \textbf{Extensibility to higher orders.} Switching from second to fourth or sixth order requires changing a single parameter in the MOLE library. The operator structure and update code remain identical.

  \item \textbf{Modularity.} Material properties and absorbing boundary conditions (e.g., UPML) are incorporated through element-wise coefficients, without modifying the operators. The same operators work on Cartesian and curvilinear grids~\cite{Boada2025} and have been applied to anisotropic elliptic problems~\cite{Boada2020}.
\end{enumerate}

We note that the one-dimensional formulation ($E_x$, $H_y$) and the two-dimensional TM$_z$ formulation ($E_z$, $B_x$, $B_y$) presented here require only the mimetic divergence and gradient operators, since the curl reduces to scalar partial derivatives in these settings. A fully three-dimensional formulation, however, requires the mimetic curl operator---also available in the MOLE library---to discretize Faraday's and Amp\`ere's laws in their full vector form.

Future work includes a systematic convergence study comparing the mimetic approach with standard FDTD across different orders of accuracy, as well as extensions to three-dimensional and curvilinear geometries.

The mimetic framework, implemented in the open-source MOLE library~\cite{Corbino2024}, provides a principled and practical alternative to FDTD for computational electromagnetics. All examples presented in this paper, along with additional electromagnetic simulations, are available in the MOLE repository.

% ══════════════════════════════════════════════════════════════

\end{document}